 \documentclass[11pt, a4paper]{article}       

\usepackage{amsmath, amsthm, amssymb,amsfonts}
 \usepackage{wasysym}
\usepackage[utf8]{inputenc}
\usepackage[T1]{fontenc}
\usepackage{url}
\usepackage{verbatim}
\usepackage{indentfirst}
\usepackage{xypic}
\usepackage{titlesec}
\usepackage{fancyhdr}
\usepackage{MnSymbol}
\usepackage{textcomp}
\usepackage{braket}
\usepackage{enumerate}
\usepackage{stmaryrd}
\usepackage{xcolor}
\usepackage{hyperref}
\usepackage{pifont}
\usepackage{float}
\usepackage{graphicx}
\usepackage[english,francais]{babel}

\usepackage{geometry}
\theoremstyle{plain}

  \newtheorem*{theorem*}{Theorem}

  \newtheorem*{proposition*}{Proposition}

\theoremstyle{definition}

\DeclareFontFamily{OT1}{pzc}{}
\DeclareFontShape{OT1}{pzc}{m}{it}{<-> s * [1.1] pzcmi7t}{}
\DeclareMathAlphabet{\mathpzc}{OT1}{pzc}{m}{it}

\renewcommand{\u}{_}
\newcommand{\m}{-}
\newcommand{\R}{\mathbb{R}}
\newcommand{\C}{\mathbb{C}}

\newcommand{\N}{\mathbb{N}}

\newcommand{\norme}[1]{\left\Vert #1\right\Vert}

\title{Un effet de moiré \\ sur les espaces symétriques de type non\m compact}

\author{Alexandre Afgoustidis\footnote{CEREMADE (UMR CNRS no. 7534), Universit\'e Paris-Dauphine, Place du Mar\'echal de Lattre de Tassigny, 75775 Paris Cedex 16, France, \emph{and} Institut de Math\'ematiques de Jussieu Paris Rive Gauche, Universit\'e Paris Diderot, B\^atiment Sophie Germain, 75013, Paris, France. Email: \texttt{afgoustidis@ceremade.dauphine.fr}}}

\date{Juin 2014}

\begin{document}
\maketitle
\begin{abstract}
Nous montrons que sur un espace symétrique $X$ de type non compact, de même que l'interférence constructive d'ondes d'Helgason dont les directions de propagation balaient le "bord" de $X$ fournit ses fonctions sphériques élémentaires, une onde d'Helgason peut être obtenue par l'interférence constructive de fonctions sphériques élémentaires dont les centres balaient un horocycle de $X$. 
\end{abstract}
%
%
%
\section{Introduction}


Lorsqu'on obtient un motif reconnaissable en superposant deux ou plusieurs figures du plan qui ne diffèrent que de petits mouvements euclidiens, on parle d'\emph{effet de moiré}. Pour des exemples venus du traitement d'images, voir \cite{Harthong, Amidror}. Indiquons un exemple d'utilisation possible en neurosciences \cite{Hubel, HW, Hub}, qui est notre motivation. 

Le long de la voie qui mène de la rétine au cortex visuel, on trouve des neurones (dans un noyau du thalamus nommé corps genouillé latéral, LGN en abrégé) dont il est raisonnable de décrire la réponse à l'aide une longueur d'onde $\lambda$, d'un point $x\u{0}$ du plan euclidien et de la fonction de Bessel sphérique $J\u{\lambda, x\u{0}}$ correspondante\footnote{$J\u{\lambda, x\u{0}}$ est la fonction $\vec{x} \mapsto \int\u{\mathbb{S}^1} \exp\left( \frac{2\pi}{\lambda} \vec{u} \cdot \vec{(x \m x\u{0})} \right) du$, définie sur le plan euclidien.} :  en envisageant l'image vue à l'instant $t$ comme une fonction à valeurs réelles $I\u{t}$ sur le plan euclidien, en négligeant certains effets de localisation en simplifiant l'intégration temporelle, le produit scalaire $\mathbb{L}^2$ entre $J\u{\lambda, x\u{0}}$ et $I\u{t}$ fournit une bonne description  pour l'activité du neurone à l'instant $t+\delta t$, où $\delta t$ est le délai de réponse.

Le LGN apporte l'information visuelle au cortex visuel primaire (V1), mais dans ce dernier la réponse des neurones est (toujours en négligeant des effets de localisation) bien mieux décrite à l'aide d'un (co\m)vecteur $\vec{p}$, en formant le produit scalaire $\mathbb{L}^2$ entre $I\u{t}$ et l' onde plane $W\u{\vec{p}} = \vec{x} \mapsto e^{i \vec{p} \cdot \vec{x}}$. Pour obtenir une explication raisonnable à la transformation des profils récepteurs entre LGN et V1, on peut s'appuyer \cite{HW,Hubel} sur l'observation suivante~: si $(x\u{i})\u{i=1..n}$ est une famille de points alignés et voisins, la somme $ \frac{1}{n} \sum J\u{\frac{2\pi}{\norme{p}}, x\u{i}}$ n'est pas sans ressembler à $W\u{\vec{p}}$. Le fait que l'arbre dendritique d'une cellule dont $W\u{\vec{p}}$ décrit bien le profil récepteur, reçoive généralement des contacts synaptiques des cellules ganglionnaires dont les profils ont à voir avec les $J\u{2\pi/\norme{\vec{p}}, x\u{i}}$, est à cet égard très spectaculaire.

Il est tentant de regarder ici le lien entre la fonction de Bessel sphérique et les ondes planes avec les lunettes de l'analyse harmonique non\m commutative, en faisant intervenir le groupe des déplacements du plan. Il est également tentant de généraliser ce résultat à d'autres espaces homogènes, par exemple en montrant que la reconstruction d'une "onde plane"  à partir de fonctions sphériques reste valable dans un espace symétrique de type non\m compact.\\ 

Cette courte note a vocation à être incluse dans ma thèse de doctorat de l'université Paris 7. Merci à mon directeur de thèse Daniel Bennequin : ses remarques sur le cerveau sont à l'origine de ce travail.

\section{Notations}

Soit $G$ un groupe de Lie semi-simple réel non compact ; nous supposerons que le centre de $G$ est fini. Adoptons les notations usuelles \cite{H1} : 

$\bullet$ Notons $\mathfrak{g}$ l'algèbre de Lie de $G$, $K$ un sous\m groupe compact maximal, $A$ et $N$ des sous\m groupes (l'algèbre de Lie de $A$ est notée $\mathfrak{a}$) fournissant une décomposition d'Iwasawa $G = KAN$ de $G$. Le sous\m groupe $N$ vient avec un choix de système de racines positives pour la paire $(\mathfrak{g},\mathfrak{a})$ ; je noterai $\rho$ la demi\m somme des racines positives associée, qui est un élément de $\mathfrak{a}^\star$.

$\bullet$ Appelons $M$ le centralisateur de $A$ dans $K$, et $B$ le quotient compact $K/M$.  

$\bullet$ Supposons que des mesures $G$\m invariantes sont fixées sur $G$ et chacun de ses sous\m groupes et quotients de façon cohérente (voir Helgason \cite{H2}), et de façon à ce que le volume total de $B$ soit $1$. Les intégrations à venir seront, sauf précision, relatives à ces mesures invariantes.\\ 

Notons $X$ l'espace symétrique de type non-compact $G/K$. Suivant Helgason \cite{H3}, nous appelons \emph{horocycle}\footnote{Poincaré disait \emph{horisphère}} une orbite dans $X$ d'un conjugué de $N$.  Si $g = k_0a_0n_0$ est un élément de $G$, le sous-groupe $gNg^{-1} = k_0 N k_0^{-1}$ ne dépend que de l'image de $K$ dans le quotient $B = K/M$ ; appelons \emph{direction} d'une orbite de $gNg^{-1}$ l'image de $k_0$ dans $K/M$. Pour $x \in X$ et $b \in B$, notons $\xi(b,x)$ l'unique horocycle de direction $b$ passant par $x$, orbite de $x$ sous le conjugué de $N$ correspondant à $b$. 

Soit $b \in B$ et $x = n a \in X$ ; posons $\Delta(x ; b) = \mathcal{A}(b^{-1} \tilde{x}) \in \mathfrak{a}$ où $\tilde{x}$ est un relevé quelconque de $x$ dans $G$ et où $\mathcal{A} : G \mapsto \mathfrak{a}$ désigne la projection $n a k \mapsto log_A(a)$ issue de la décomposition d'Iwasawa. La quantité $\Delta(x;b)$ ne dépend que de l'horocycle $\xi(x,b)$. 

Lorsque $G$ est le groupe $SU(1,1)$ et agit par homographies sur le disque unité ouvert $\mathbb{D}$ de $\C$, le stabilisateur de $0$ est un compact maximal $K$, isomorphe à $SO(2)$ ; les horocycles de $\mathbb{D}$ (qui ne dépendent que du choix de $K$) sont les cercles tangents au bord du disque ; il est naturel d'identifier la direction d'un horocycle avec le point de tangence. 

\begin{figure}[here]
\begin{center}
\includegraphics[width=0.34\textwidth]{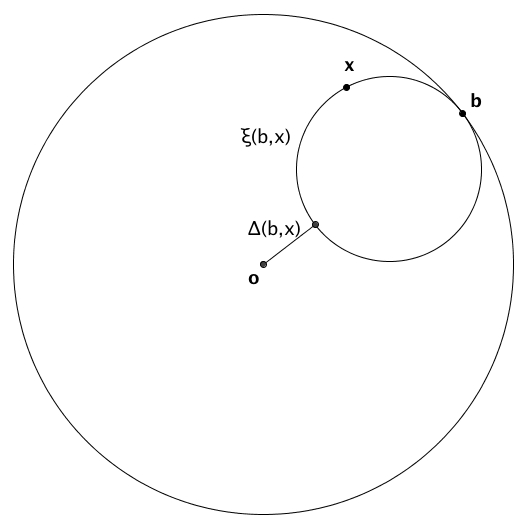}
\caption{Un horocycle de $\mathbb{D}$. }
\end{center}
\end{figure}
\newpage
\section{Transformée de Fourier-Helgason}

Pour $\lambda \in \mathfrak{a}^\star$ et $b \in B$, définissons
\begin{align*} e_{\lambda,b} : X \rightarrow & \R \\
x \mapsto & e^{\braket{i\lambda + \rho \ | \ \Delta(b;x)}}. \end{align*}
La fonction $e_{\lambda,b}$ est constante sur la famille des horocycles de direction $b$ et joue pour l'analyse harmonique $G-$invariante sur $X$ le rôle que jouent les ondes planes pour l'analyse de Fourier sur l'espace euclidien. 

\begin{figure}[here]
\begin{center}
\includegraphics[width=0.4\textwidth]{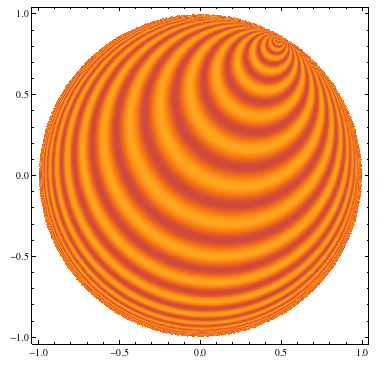}
\caption{Cette figure représente les \emph{lignes de phase} d'une onde d'Helgason (mais ne tient pas compte de sa croissance à l'infini).}
\end{center}
\end{figure}

Si $f$ est une fonction de $X$ dans $\C$, la transformée de Fourier-Helgason de $f$ est la fonction 
\[ \hat{f} : (\lambda \in \mathfrak{a}^\star, b \in K/M) \mapsto \int \u{X} e\u{\m \lambda,b}(x) f(x) dx \]
définie sur le sous-ensemble de $\mathfrak{a}^\star \times K/M$ sur lequel cette intégrale converge, par exemple sur $\mathfrak{a}^\star \times K/M$ tout entier si $f$ est une fonction lisse à support compact sur $\mathcal{X}$. Lorsque le chapeau s'avèrera trop étroit pour que la notation soit lisible,j'écrirai $\mathcal{F}(f)$ pour $\hat{f}$. 

Lorsque $f$ est une fonction intégrable sur $X$, il n'est plus évident que cette intégrale converge pour tout $(\lambda,b) \in \mathfrak{a}^\star \times K/M$ ; cependant, on peut montrer (\cite{H3}, p. 209) qu'elle est définie sur $\mathfrak{a}^\star \times B\u{0}$, où $B\u{0}$ est un sous-ensemble de pleine mesure de $B$.

Il est utile pour ce qui suivra d'étendre la définition de la transformation de Fourier aux distributions sur l'espace $X$, et de disposer d'analogues non\m euclidiens de l'espace de Schwartz $\mathcal{S}(X)$ et de l'espace  $\mathcal{S}'(X)$ des distributions tempérées ; les définitions correspondantes (\cite{H3}, pp. 214 et suivantes) ont été renvoyées en fin de texte. Mais même dans le cas où $f$ n'est qu'intégrable, lorsque $\hat{f}$ est intégrable relativement à la mesure $\left( |\mathbf{c}(\lambda)|^{\m 2} d\lambda \right) \otimes db$ sur  $\mathfrak{a}^\star \times K/M$ (cette mesure fait intervenir la fonction $\mathbf{c}$ d'Harish\m Chandra), on dispose de la formule d'inversion suivante, valable pour presque tout $x \in X$ : 
\[ f(x) = \frac{1}{w} \int\u{\mathfrak{a^\star} \times B} \hat{f}(\lambda,b)\ e\u{\lambda,b}(x) \  \mathbf{c}(\lambda)^{\m 2} d\lambda  db \]
où $w$ est l'ordre du groupe de Weyl $W(\mathfrak{g}, \mathfrak{a})$ ; la fonction $\mathbf{c}$ est analytique sur le complémentaire dans $\mathfrak{a}^\star$ d'une réunion finie d'hyperplans, mais je n'aurai pas besoin de détailler ici. Et les questions sur la transformée de Fourier euclidienne trouvent souvent des termes pour être posées sur la transformation de Fourier-Helgason de $X$ (formule de Plancherel, théorème de Paley-Wiener...) ; les réponses à ces questions rappellent souvent le cas euclidien, bien qu'il y ait des différences importantes dues à la courbure de $X$ (et à la croissance à l'infini des ondes d'Helgason qu'elle impose).

\section{Fonctions sphériques élémentaires}
Les ondes $e\u{\lambda,b}$ donnent un air familier à l'analyse harmonique $G$\m invariante sur $G/K$, peinte en détail dans les ouvrages d'Helgason ; avant qu'Helgason ne révèle leur rôle, le fait que la fonction obtenue par l'interférence constructive de toutes les ondes de longueur d'onde $\lambda$ soit la fonction sphérique $\varphi_\lambda$ avait été crucial dans le programme d'Harish-Chandra pour l'étude du dual tempéré de $G$ \cite{HC}.

Pour tout $\lambda \in \mathfrak{a}^\star$, la fonction  
\begin{align*} \varphi\u{\lambda} : G \rightarrow & \ \C \\
g \mapsto & \int \u{B} e\u{\lambda,b}(gK) \ db \end{align*}
vaut $1$ en $o$, est $K$\m bi\m invariante et fonction propre de tous les opérateurs différentiels $G$\m invariants sur $G$ : c'est une fonction sphérique élémentaire de $G$. Les seules fonctions vérifiant les trois propriétés de la phrase précédente sont les $\varphi\u{\lambda}$, $\lambda \in \mathfrak{a}^\star$, et deux fonctions de ce type sont égales lorsque les éléments de $\mathfrak{a}^\star$ qui y apparaissent se correspondent par l'action d'un élément du groupe de Weyl. 

Nous n'en aurons pas besoin, mais rappelons que la version $K$\m invariante de la transformation de Fourier\m Helgason a mené à la définition de la fonction $\mathbf{c}$ et ouvert la marche triomphale d'Harish\m Chandra vers la formule de Plancherel pour $G$  : si $f$ est une fonction lisse à support compact et $K$\m bi\m invariante sur $G$, notons $\tilde{f}(\lambda) = \int\u{G} f(g) \varphi\u{\m \lambda}(g) dg$ pour tout $\lambda \in \mathfrak{a}^\star$, alors
\[f \mapsto \tilde{f} \ \text{s'étend en une isométrie de} \ \mathbb{L}^2(K\backslash G /K) \ \text{vers} \ \mathbb{L}^2 \left(\mathfrak{a}^\star\u{+}, |\textbf{c}(\lambda)|^{\m 2}d\lambda \right).   \]

\section{Effet de moiré}

Dans les quelques lignes qui suivent, il s'agit de montrer que les ondes planes d'Helgason sont, réciproquement, reconstruites par l'interférence entre des fonctions sphériques dont les centres sont disposés le long d'un horocycle. 

Précisément, fixons une "longueur d'onde" $\lambda \in \mathfrak{a}^\star$, un point du bord $-$ disons l'image  $b_0$  dans $B$ de l'identité de $K$, et un point $x$ de $X$. Pour tout élément $y$ de $X$, notons $\phi_\lambda^{[y]}$ l'unique membre de l'espace $\mathcal{E}_\lambda(X)$ du chapitre de \cite{H3}\footnote{Il s'agit de l'espace propre commun aux opérateurs différentiels $G$\m invariants sur $X$ qui contient $\varphi\u{\lambda}$} valant $1$ en $y$, insensible aux translations à gauche et à droite de la variable le long d'un élément du stabilisateur de $y$ . Nous allons montrer l'égalité suivante : 

\begin{equation} \label{E1} \hspace{\m 5.4cm} \text{\textbf{Résultat principal.}} \hspace{2.3cm} \int_{\xi(b_0,0)} \varphi^{[y]}_\lambda(x) \ dy = e_{\lambda,b_0}(x). \end{equation}

\begin{figure}[here]
\begin{center}
\includegraphics[width=0.4\textwidth]{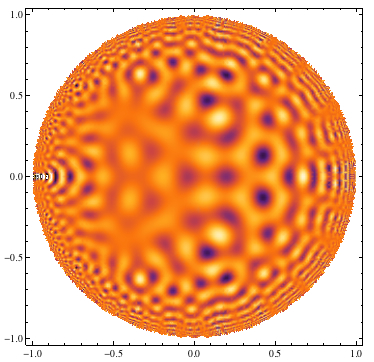} \includegraphics[width=0.4\textwidth]{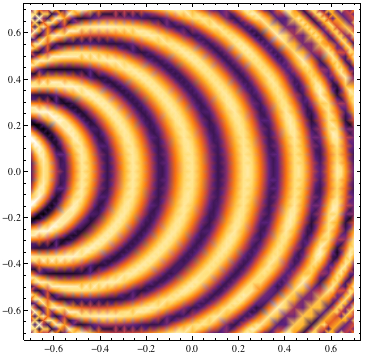}
\caption{Somme de 5, puis (détail de la somme) de 60 fonctions sphériques centrées sur l'horocycle $\xi(\m 1, 0)$ de $\mathbb{D}$. }
\end{center}
\end{figure}

Si $y \in X$ est l'image de l'origine de $G/K$ par l'élément $g_y$ de $G$, la fonction $\phi_\lambda^{[y]}$ n'est autre que $z \mapsto \phi_\lambda(g_y^{-1} \cdot z)$. Lorsque $y$ est sur l'horocycle $\xi(b_0,0)$, l'élément $g_y$ peut être choisi dans $N$ et $g_y^{-1} \cdot x$ est sur l'horocycle $\xi(b_0,x)$. Il s'agit donc de montrer l'égalité

\[ \int_{\xi(b_0,x)} \phi_\lambda = e_{\lambda,b_0}(x). \]

%

Mais la formule intégrale d'Harish-Chandra déjà mentionnée (\cite{HCPNAS}) permet d'écrire 

\begin{equation}  \label{HC2} \int_{\xi(b_0,x)} \varphi_\lambda =  \int_{\xi(b_0,x)} dy \int_B e_{\lambda,b}(y) db. \end{equation}

Echanger les sommations dans \eqref{HC2} ferait apparaître l'intégrale sur le bord $B$ de la transformée de Fourier-Helgason de l'indicatrice de l'horocycle $\xi(b_0,x)$, ce qui nous mènerait près de la conclusion souhaitée. Mais cela ferait aussi apparaître des intégrales divergentes. Notons cependant le résultat suivant.  \\

%

\textbf{Lemme :} \emph{la transformée de Fourier de la distribution de Dirac sur l'horocycle $\xi(b_0, x)$ (vue comme une distribution tempérée sur $\mathfrak{a}^\star \times B$) est :}

\begin{equation} \label{TF} \left[ \lambda \mapsto e_{\lambda,b_0}(x) \right] \otimes  \delta_{b=b_0}. \end{equation} 

Cela rappelle le fait que de la distribution de Dirac sur une droite de $\R^2$ qui passe par l'origine a pour transformée de Fourier la distribution de Dirac sur la droite orthogonale qui passe par l'origine ; notons cependant que la direction privilégiée qui apparaît dans la transformée de Fourier\m Helgason est ici la direction de l'horocycle lui\m même plutôt qu'une direction "orthogonale".\\

\emph{Preuve : } Notons $T_{b_0,x} \in \mathcal{S}'(X)$ la distribution \eqref{TF}. Nous devons vérifier que pour tout $\psi \in \mathcal{S}(X)$,
\[ \braket{T_{b_0,x} | \hat{\psi}} = \int_{\xi(b_0, x)} \psi.\] 
Mais  $ \braket{T_{b_0,x} | \hat{\psi}} = \int_{\mathfrak{a}^\star} e_{\lambda,b_0}(x) \hat{\psi}(\lambda,b_0) d\lambda$, soit bien sûr
\begin{align*}
\braket{T_{b_0,x} | \hat{\psi}} =& \int_{\mathfrak{a}^\star} e_{\lambda,b_0}(x) \left[ \int_X  e_{-\lambda,b_0}(y) \psi(y) dy \right] d\lambda \\
=& \int_{\mathfrak{a}^\star}\left[ \int_X   e_{\lambda,b_0}(x)  e_{-\lambda,b_0}(y) \psi(y) dy \right] d\lambda \\
=& \int_{\mathfrak{a}^\star}\left[ \int_X   e^{ \braket{i \lambda + \rho | (\Delta(b_0; x) -\Delta(b_0; y) }} \psi(y) dy \right] d\lambda.
\end{align*}

Notons maintenant $\sigma$ la fonction $y \mapsto \Delta(b_0 ; x) - \Delta(b_0 ; y)$ de $X$ dans $\mathfrak{a}$. Les ensembles de niveaux de $\sigma$ sont les horocycles de direction $b_0$. En munissant chaque horocycle de la mesure induite par la métrique invariante de $X$ et  en utilisant la formule de la co\m aire, on obtient
\begin{align*}
\braket{T_{b_0,x} | \hat{\psi}} =& \int_{\mathfrak{a}^\star}\left[ \int_X   e^{ \braket{i \lambda + \rho | (\Delta(b_0; x) -\Delta(b_0; y) }} \psi(y) dy \right] d\lambda \\
=&  \int_{\mathfrak{a}^\star} \int_{\mathfrak{a}}   e^{ \braket{i \lambda | u}} \left\{ \int_{\sigma^{-1}(u)} \psi e^{\braket{\rho | u}} \right\} du d\lambda. 
\end{align*}
Il ne reste plus qu'à reconnaître dans la dernière expression l'intégrale sur $\mathfrak{a}^\star$ de la transformée de Fourier euclidienne de
 \begin{align*}
\Psi :  \mathfrak{a} \rightarrow & \ \C \\
 u \mapsto & \int_{\sigma^{-1}(u)} \psi(z) e^{\braket{\rho | u}} dz_u
 \end{align*}
pour obtenir, grâce à l'inversion de Fourier euclidienne, 
\[ \braket{T_{b_0,x} | \hat{\psi}} = \Psi(0) = \int_{\sigma^{-1}(0) = \xi(b_0,x)} \psi \]
\begin{flushright} C.Q.F.D. \end{flushright}
~\\
%
\ \ \ Revenons à l'effet de moiré \eqref{E1}. L'échange d'intégrales dans  \eqref{HC2} étant manifestement souhaitable, remarquons que pour tout élément $\vartheta$ de $\mathcal{S}(X)$, $\mathcal{F} \left\{p \mapsto \vartheta(p) \ \mathbf{1}_{\xi(b_0,p)} \right\}$ est définie et mesurable sur un ensemble de la forme $\mathfrak{a}^\star \times B\u{0}$ où $B\u{0} \subset B$ est de pleine mesure, et continue en $\lambda$ sur cet ensemble. Nous pouvons donc fixer $\lambda \in \mathfrak{a}^\star$, puis utiliser le fait que t $(p, \beta) \mapsto \vartheta(p) e\u{\lambda,b}(p)$ est intégrable sur $X \times B$ pour échanger les intégrales, et écrire

\begin{align} \int_{\xi(b_0,x)}  \vartheta \varphi\u{\lambda}= \int_{\xi(b_0,x)} dy \ \vartheta(y) \int_B e_{\lambda,b}(y) \ db =& \int_B db \int_X \vartheta(p) \ \mathbf{1}_{\xi(b_0,x)}(p) \ e_{\lambda,b}(p)  dp \\
=& \int_B \mathcal{F} \left\{ \vartheta \ \mathbf{1}_{\xi(b_0,\cdot)} \right\} (\lambda,b) db. \end{align}


Dans la dernière égalité la quantité $\int_B \mathcal{F} \left\{ \vartheta \ \mathbf{1}_{\xi(b_0,x)} \right\} (\lambda,b) db$ définit une fonction continue de $\lambda$ : en effet, la fonction (à valeurs complexes) $\phi\u{\lambda}$ est de module inférieur ou égal à $1$ pour $\lambda \in \mathfrak{a}^\star$ et donc

\begin{align*}
 \int_{B}  \left| \mathcal{F} \left(\vartheta \ \mathbf{1}_{\xi(b_0,x)}\right) \right|(\lambda,b) db  \leq & \int_{X} \left| \vartheta \ \mathbf{1}_{\xi(b_0,x)} \right| \left(\int_{B} e^{\braket{\rho \ | A(x,b)}} db \right) dx \\
 \leq &  \int_{\xi(b_0,x)} \left| \vartheta \right|,\\
\end{align*}

ce qui permet de déduire la continuité annoncée des théorèmes ordinaires d'intégration et du fait que pour $x$ fixé dans $X$, l'application $\lambda \mapsto \varphi\u{\lambda}(x)$ est continue. 

Faisons maintenant tendre $\vartheta$ vers la fonction constante égale à $1$ sur $X$, au sens des distributions tempérées ; les propriétés de continuité (A1) et (A2)  rappellées en fin de texte permettent de s'appuyer sur le calcul précédent pour conclure qu'à $\lambda$ fixé, le nombre complexe $\int_B \mathcal{F} \left\{ \vartheta \ \mathbf{1}_{\xi(b_0,x)} \right\} (\lambda,b) db$ tend vers $e_{\lambda,b_0}(x)$. Nous avons démontré l'égalité \eqref{E1}, ce qui était le but de cette courte note.

\section{Distributions tempérées sur un espace symétrique}

\subsection{Classe de Schwartz de X.}

Soit $\mathbf{D}(G)$ l'algèbre des opérateurs différentiels sur $G$ invariants à gauche par les translations de $G$, et $\bar{\mathbf{D}}(G)$ l'algèbre des opérateurs invariants à droite.

Rappelons que tout élément de $G$ peut s'écrire comme un produit $k\u1 a k\u 2$ où $k\u 1$ et $k\u 2$ appartiennent à $K$ et $a$ à $A$, et que deux décompositions de ce type d'un même élément voient leurs éléments $a$ différer de l'action d'un élément du groupe de Weyl. Poser $|g| = |\log(a)|$ (le membre de droite fait référence à une norme euclidienne sur $\mathfrak{a}$) conduit à la définition suivante : une fonction lisse $f$ sur $G$ est à décroissance rapide si pour tous $\ell \in \N$,  $L \in \mathbf{D}(G)$ et $R \in \bar{\mathbf{D}}(G)$, 
\begin{equation}\label{Swa} \text{sup}_{g \in G} \left| (1+|g|)^{\ell} \Xi(g)^{-1} (LR f)(g) \right| \end{equation}
est fini. 

Dans \eqref{Swa}, la fonction $\Xi$ est la fonction sphérique  $\varphi\u 0$ ; elle vérifie [HC, théorème 3] : 

\[ \Xi(g) \leq c (1+|g|)^d e^{\m \braket{\rho | \log a }} \]

où $c$ est un réel positif et $d$ un entier naturel. 

Ensemble, les fonctions à décroissance rapide sur $G$ forment l'espace de Schwartz $\mathcal{S}(G)$~; celles qui sont invariantes à droite par $K$ forment l'espace $\mathcal{S}(X)$. Les grandeurs \eqref{Swa} permettent de munir ces espaces de topologies d'espaces de Fréchet ; je noterai $\mathcal{S}'(X)$ le dual de $\mathcal{S}(X)$, espace des \emph{distributions tempérées} sur $X$.

\subsection{Classe de Schwartz de $\mathfrak{a}^\star \times K/M$ et continuité de la transformée de Fourier.}
 
L'image de $\mathcal{S}(X)$ par la transformée de \m Helgason Fourier est formée de fonctions lisses sur $\mathfrak{a}^\star \times K/M$ vérifiant (\cite{H3}, chap. 3, thm 1.10) : 

\begin{equation}\label{Swa2} \text{For each} \ P \in \R[X,Y] \ \text{and every}\ \ell \in \N,  \ \ \underset{\lambda,b}{\text{sup}} \left| (1+|\lambda|)^\ell \left(P(\Delta_{K/M},\Delta_{\mathfrak{a}^\star})\cdot g\right)(\lambda,b)  \right| < \infty \end{equation}

où $\Delta\u{K/M}$ et $\Delta\u{\mathfrak{a}^\star}$ sont les opérateurs de Laplace\m Beltrami sur $B$ et $\mathfrak{a}^\star$. 

Notons $\mathcal{S}(\mathfrak{a}^\star \times K/M)$ l'ensemble des fonctions lisses sur $\mathfrak{a}^\star \times B$ satisfaisant \eqref{Swa2} ; il vient naturellement avec des semi\m normes qui le munissent d'une topologie d'espace de Fréchet~; la transformation de Fourier\m Helgason est bien sûr continue et injective de $\mathcal{S}(X)$ dans $\mathcal{S}(\mathfrak{a}^\star \times K/M)$. 

Harish\m Chandra, puis Helgason ont montré qu'elle induit un homéomorphisme entre les sous\m espaces de $\mathcal{S}(X)$ et $\mathcal{S}(\mathfrak{a}^\star \times K/M)$ qui rassemblent les éléments $K\m$invariants (\cite{H3}, th. 1.17 ; voir aussi Anker \cite{Anker}). Eguchi \cite{Eg, EO} a vérifié (avec Okamoto) que c'est un homéomorphisme. 

Définir les distributions tempérées sur $\mathfrak{a}^\star \times K/M$ et leur transformée de Fourier de la manière habituelle\footnote{L'espace $\mathcal{S}'(\mathfrak{a}^\star \times K/M)$ des distributions tempérées est le dual topologique $\mathcal{S}(X)$. Si $T$ une distribution tempérée sur $X$, $\hat{T}$ est la distribution tempérée $\psi = \hat{\varphi} \in \mathcal{S}(\mathfrak{a}^\star \times B) \mapsto \braket{ T \ | \ \varphi }$ sur $\mathfrak{a}^\star \times B$. } donne bien sûr : 
\begin{equation}
T \mapsto \hat{T}  \ \text{est un homéomorphisme entre} \ \mathcal{S}'(X) \ \text{et} \ \mathcal{S}'(\mathfrak{a}^\star \times K/M). \tag{A1}
\end{equation}
Nous en avons eu besoin pour démontrer notre résultat ; la remarque suivante est aussi nécessaire : si $T$ est une distribution tempérée sur $\mathfrak{a}^\star \times K/M$, on définit une distribution $U$ sur $\mathfrak{a}^\star$ ("intégrale de $T$ sur $K/M$")  en posant, pour tout $\zeta \in \mathcal{S}(\mathfrak{a}^\star)$, $\braket{U \ | \ \zeta} = \braket{ T \ | \ \zeta \otimes 1\u{B} }$. Et bien sûr,
\begin{equation}
\text{L'application} \ T \mapsto U  \ \text{est continue de} \ \mathcal{S}'(\mathfrak{a}^\star \times K/M) \ \text{vers} \ \mathcal{S}'(\mathfrak{a}^\star).\tag{A2}
\end{equation} 



\begin{thebibliography}{10}


%

\bibitem{Amidror}

{I. Amidror}, \emph{The theory of the moiré phenomenon}, Computational Imaging and Vision  (2000), vol. 15.

\bibitem{Anker}
{J-Ph. Anker} , \emph{The spherical Fourier transform of rapidly decreasing functions}, J. Funct. Anal. \textbf{96-2} (1991), pp. 331–349.

\bibitem{Eg}
{M. Eguchi}, \emph{Asymptotics of Eisenstein Integrals} Journal of Functional Analysis \textbf{34} (1979), pp. 167-216 

\bibitem{EO}

 {M. Eguchi \& K. Okamoto}, \emph{The Fourier Transform of the Schwartz Space on a Symmetric Space}, Proc. Japan Acad. Ser. A \textbf{53 -7}  (1977), pp. 237-241

\bibitem{HCPNAS}

{Harish-Chandra}, \emph{Spherical functions on a semisimple Lie group}, Proc Natl Acad Sci U S A. \textbf{43(5)} (1957), pp. 408–409.

\bibitem{HC}
{Harish-Chandra}, \emph{Spherical functions on a semisimple Lie group. I}, American Journal of Mathematics \textbf{80} (1958), pp. 241–310

\bibitem{Harthong}

{J. Harthong}, \emph{Le moiré}. Publications de l'IRMA, Strasbourg, vol. 93 (1981).

\bibitem{H1}
{S. Helgason}, \emph{Differential geometry, Lie groups, and symmetric spaces}. Academic press, vol. 80 (1979).

\bibitem{H2}
{S. Helgason}, \emph{Groups \& Geometric Analysis: Radon Transforms, Invariant Differential Operators and Spherical Functions}. Academic press (1984).

\bibitem{H3}
{S. Helgason}, \emph{Geometric analysis on symmetric spaces}, 2nd Ed. American Mathematical Society, vol. 39 (2008).

\bibitem{Hubel}

{D. Hubel}, \emph{Eye, Brain and Vision}. Scientific American Library (1995).

\bibitem{HW}

{D. H. Hubel and T. N. Wiesel}, \emph{Receptive fields, binocular interaction and functional architecture in the cat's visual cortex}, Journal of Physiology \textbf{160-1} (1962), pp. 106-154.

\bibitem{Hub}

{D. Hubel}, \emph{Transformation of information in the cat's visual system}, Proceedings of the International Union of Physiological Sciences, XXII International Congress, Leiden (1962).


\end{thebibliography}
\end{document}